# Polyhedral Spline Finite Elements

**Zhao Guohui**

**School of Mathematical Sciences, Dalian University of Technology**

**Dalian City, Liaoning Province, P.R. China, 116024**

**Zhao Guohui, School of Mathematical Sciences, Dalian University of Technology,**

**No.2 Linggong Road, Ganjingzi District, Dalian City, Liaoning Province,**

**P.R. China, 116024, ghzhao@dlut.edu.cn, 86-0411-84708351-8315**

# Polyhedral spline finite elements

Zhao Guohui

School of Mathematical Sciences, Dalian University of Technology,

No. 2 Linggong Road, Ganjingzi District, 116024 Dalian Liaoning, P.R. China

**Abstract**

Spline functions have long been used in numerically solving differential equations. Recently it revives as isogeometric analysis, which uses NURBS for both parametrization and element functions. In this paper, we introduce some multivariate spline finite elements on general partitions. These splines are constructed by mollifying splines of low smooth order with B-splines. They have both high smooth order and strong approximation power. The low smooth order splines can be chosen flexibly, for example polynomials on individual cells of general polyhedral partitions. The obtained spline elements are suitable for adaptive computation.

# 1. Polyhedral splines: box splines, simplex splines and cone splines

This paper is mainly concerned with an application of multivariate splines. It is a generalization of isogeometric analysis. So some basic results of multivariate splines is listed below. For detailed introduction, please refer to the references and the internet. Multivariate spline theory has been studied extensively. It is used in approximation, CAGD, FEM, etc. For the sake of simplicity, we restrict to the bivariate case. Given a polygonal domain $\Omega$, we partition it with irreducible algebraic curves into finite cells. The partition is denoted by $\Delta$. A function $f(x)$ defined on $\Omega$ is a spline of degree $n$ and smooth order $r$ if $f \in c^r(\Omega)$ and on each cell $f(x)$ is a polynomial of degree $n$. The set of all splines of degree $n$ and smooth order $r$ is a linear space and is denoted by $s_n^r(\Omega, \Delta)$.

In this paper we mainly consider polyhedral splines, which are defined on partitions formed by line segments, that is, cells of the partition are polygons. The well-known polyhedral splines are box splines, simplex splines and cone splines. We introduce them by homogenization and symmetrization.

A univariate B-spline is defined by $B(x) = (x_{n+1} - x_0)[x_0, x_1, \cdots, x_{n+1}](\cdot - x)_+^n$. Homogenize $B(x)$ by replacing $x$ and $x_i$ with $\dfrac{x}{y}$ and $\dfrac{x_i}{y_i}$, respectively and then simplifying, we obtain cone splines

$$C(x,y) = \sum_{i=0}^{n+1} \frac{\begin{vmatrix} x, y \\ x_i, y_i \end{vmatrix}_*^n}{\prod_{j \neq i} \begin{vmatrix} x_j, y_j \\ x_i, y_i \end{vmatrix}}, \quad \text{or} \quad C(x,y) = \sum_{i=0}^{n+1} \frac{\begin{vmatrix} x, y, 1 \\ x_i, y_i, 1 \\ x_k, y_k, 1 \end{vmatrix}_*^n}{\prod_{j \neq i} \begin{vmatrix} x_j, y_j, 1 \\ x_i, y_i, 1 \\ x_k, y_k, 1 \end{vmatrix}}$$

Simplex splines are obtained by symmetric combination of cone splines:

$$\sum_{j \neq i} (\frac{\begin{vmatrix} x,y,1 \\ x_i,y_i,1 \\ x_i,y_i,1 \end{vmatrix}_*^n}{\prod_{j \neq i} \begin{vmatrix} x_k,y_k,1 \\ x_j,y_j,1 \\ x_i,y_i,1 \end{vmatrix}} + \sum_{\begin{vmatrix} x_k,y_k,1 \\ x_j,y_j,1 \\ x_l,y_l,1 \end{vmatrix}>0} \frac{\begin{vmatrix} x,y,1 \\ x_j,y_j,1 \\ x_k,y_k,1 \end{vmatrix}_*^n}{\prod_{\substack{l \neq j \\ l \neq k}} \begin{vmatrix} x_l,y_l,1 \\ x_j,y_j,1 \\ x_k,y_k,1 \end{vmatrix}})$$

Box splines can be obtained similarly. For $n \geq d$ nonzero direction vectors $v_k \in Z^d \setminus 0$, collected in a matrix form $E = [v_1 \cdots v_n]$, a d-variate box spline is defined recursively as follows:

$$M_E(x) = \begin{cases} 1/|\det(E)|, & x = \sum_{k=1}^d t_k v_k, t_k \in [0,1) \\ 0, & \text{otherwise} \end{cases}, \text{ if } d = n$$

and

$$M_E(x) = \int_0^1 M_{E \setminus v_n}(x - t v_n) dt, \quad n > d.$$

Box splines can be evaluated recursively:

$$(n-d)M_E(x) = \sum_{k=1}^n t_k M_{E \setminus v_k}(x) + (1-t_k) M_{E \setminus v_k}(x - v_k),$$

where $x = \sum_{k=1}^n t_k v_k$.

and the directional derivative $D_{v_k} M_E(x) = M_{E \setminus v_k}(x) - M_{E \setminus v_k}(x - v_k)$.

Theorem  Let $M_E$ and $M_F$ be two box splines. Then

$$\int M_E(x) M_F(x+y) dx = M_{E \cup F}(2m_E + y)$$

where $m_E = (v_1 + \cdots + v_n)/2$.

Theorem  Let $M_E$ and $M_F$ be two box splines. Let $X \subset E$ and $Y \subset F$ be subsets of r and s elements, respectively. Then

$$\int D_X M_E(x) D_Y M_F(x+y) dx = (-1)^r D_{X \cup Y} M_{E \cup F}(2m_E + y).$$

The relation among these splines are just the relation among box, simplex and cone, which holds for general polyhedra, and also for general polyhedral splines, that is, polyhedral splines can be expressed as linear combination of cone splines.

## 2. Polyhedral spline finite elements

Given a polygonal domain $\Omega$ with partition $\Delta$ formed by line segments. The cells of the partition are all polygons. Centroidal Voronoi tessellation is usually used to generate such partitions. Let $f_{jk}^i = x^j y^k, j+k \leq n$ on cell $\Delta_i$ and $f_{jk}^i = 0$ otherwise. Select a B-spline $B(x,y)$, such as a box spline. Let $B_{jk}^i = B * f_{jk}^i$. The linear space formed by linear combination of $B_{jk}^i$ has the following properties:

1. $B_{jk}^i$ is a B-spline;
2. The space contains all polynomials of degree at most n;
3. $B_{jk}^i$ is at least as smooth as B;
4. $B_{jk}^i$ is locally supported;
5. $B_{jk}^i$ is a linear combination of cone splines;
6. Quasi-interpolation can be constructed.

The B-spline can be chosen freely. The tensor product of univariate B-splines is a good choice. We call $B_{jk}^i$ polyhedral spline finite elements. Usually these B-splines $B_{jk}^i$ are linearly independent, but they may be linearly dependent. For example, the B-splines in $s_2^1(\Delta_2)$ are linearly dependent. Whether these B-splines are linearly independent or not, quasi-interpolation operators can be constructed for them. As $f_{jk}^i = x^j y^k$, $j+k \leq n$, we may start with some Lagrange interpolation operator $L(g) = \sum \lambda_{jk}^i(g) f_{jk}^i$, then convolute with the selected B-spline, resulting in an injective mapping of polynomials generally. This injective mapping can be turned into a quasi-interpolation operator, as is done in spline quasi-interpolation.

$f_{jk}^i$ may be negative in its domain. We can select new base polynomials by translation so that they are positive in their respective domains. The corresponding B-splines $B_{jk}^i$ are nonnegative. $f_{jk}^i$ can be piecewise polynomials, too if only they have good approximation property.

The support of $B_{jk}^i$ is Minkowski sum of supports of $f_{jk}^i$ and the mollifier. The partition and corresponding piecewise expressions of $B_{jk}^i$ can be obtained easily, which will be disscussed in the next section. Usually simplex splines mean small support.

The partition of $\Omega$ may be given flexibly, for example, the partition in the following figure is ok:

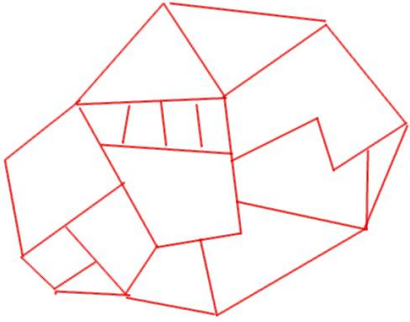

**3. Evaluation, differentiation, integration and inner product**

Polyhedral spline finite elements are B-splines. They have explicit expressions. Algorithms for exact evaluation, differentiation, integration and inner product can be given to these elements. The inner product of polyhedral spline finite elements with general functions can be computed numerically as usual.

According to the spline theory, convolution of cone splines is still a cone spline, and piecewise expressions of the convolution can be obtained easily. $f_{jk}^i$ and B-spline $B(x,y)$ are linear combinations of cone splines, so $B_{jk}^i$ is a linear combination of cone splines and has explicit piecewise expressions.

For a given direction $e = (e_1, e_2, \cdots, e_m)$, the direction integral operator is defined as usual: $J_e f(x) = \int_{-\infty}^{0} f(x+te)dt$. $J_e$ is the inverse derivative along the direction $e$.
Let $J_k = J_{e_k}, k = 1, \cdots, m$ where $e_k$ are the unit coordinate axis directions.
$$I_e f(x) = \int_{-1}^{0} f(x+te)dt = J_e f(x) - J_e f(x-e) = \nabla J_e f,$$
$I_e$ is the difference of the inverse derivative along the direction $e$.
Let $I_x = I_{e_1}, I_y = I_{e_2}, I_z = I_{e_3}$, and $f_{jk}^i = x^j y^k$, $j+k \le n$ on cell $\Delta_i$ and $f_{jk}^i = 0$ otherwise. Let $B_{jk}^i = I_x^{a+2} I_y^{b+2} f_{jk}^i$. Then $B_{jk}^i$ is a B-spline of degree $n+a+b+4$, which is a linear combination of cone splines differentiable up to (a,b) order, and
$$\frac{\partial B_{jk}^i}{\partial x}(x,y) = I_x^{a+1} I_y^{b+2} f_{jk}^i(x,y) - I_x^{a+1} I_y^{b+2} f_{jk}^i(x-1,y)$$

$$\frac{\partial B^i_{jk}}{\partial y}(x,y) = I_x^{a+2} I_y^{b+1} f^i_{jk}(x,y) - I_x^{a+2} I_y^{b+1} f^i_{jk}(x,y\text{-}1)$$

For example, the quadrilateral ABCD = cone FAE - cone FBE - cone FDE + cone FCE or $\chi_{ABCD} = \chi_{FAE} - \chi_{FBE} - \chi_{FDE} + \chi_{FCE}$, where $\chi$ is the characteristic function of the corresponding set. So $f^i_{jk} x^j y^k \chi_{ABCD} = x^j y^k \chi_{FAE} - x^j y^k \chi_{FBE} - x^j y^k \chi_{FDE} + x^j y^k \chi_{FCE}$.

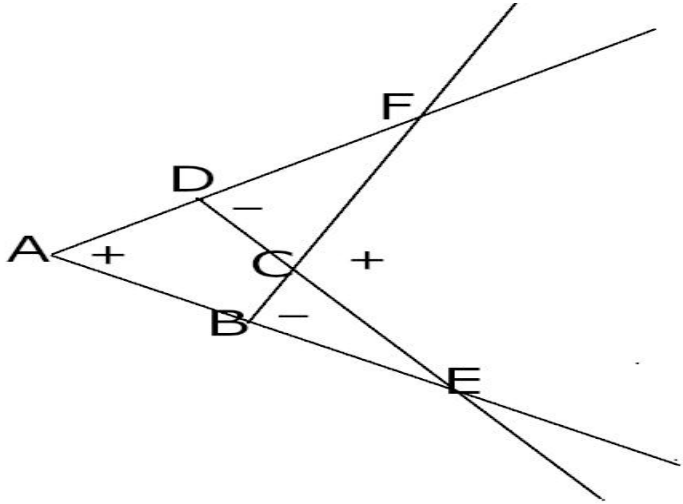

$f^i_{jk}$ is a linear combination of cone splines. To get the explicit expression of $B^i_{jk}$, it is enough to consider cone splines.

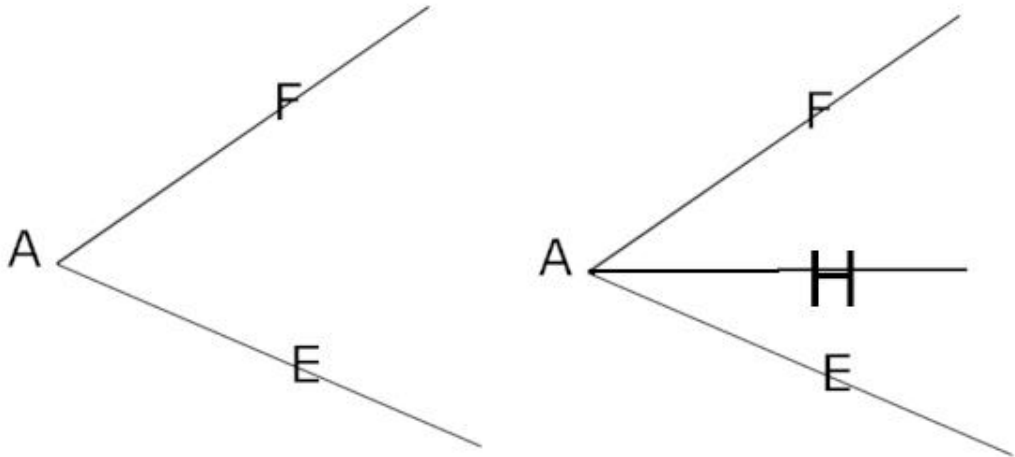

Suppose $f = x^j y^k$ on the cone FAE above and $f = 0$ otherwise. Then $I_x f$ is also a cone spline with two pieces on cones FAH and EAH where AH is the horizontal line. If the equations of AF and AE are $x + ay + b = 0$ and $x + cy + d = 0$, then

$$J_x f(x,y) = \int_{-ay-b}^{x} t^j y^k \, dt = \frac{1}{j+1}(x^{j+1} - (-ay-b)^{j+1}) \, y^k \quad \text{on cone FAH}$$

$$J_x f(x,y) = \int_{-cy-d}^{x} t^j y^k \, dt = \frac{1}{j+1}(x^{j+1} - (-cy-d)^{j+1}) \, y^k \text{ on cone EAH}$$

$J_y f$ can be obtained similarly. $B_{jk}^i$ can be obtained inductively.

The above introduction can be turned into an algorithm:
  Algorithm for evaluation, differentiation, integration and inner product:
  1 Turn a polygon into combination of cones, e.g. the pentagon in the following figure is 1-2-3-4+5.

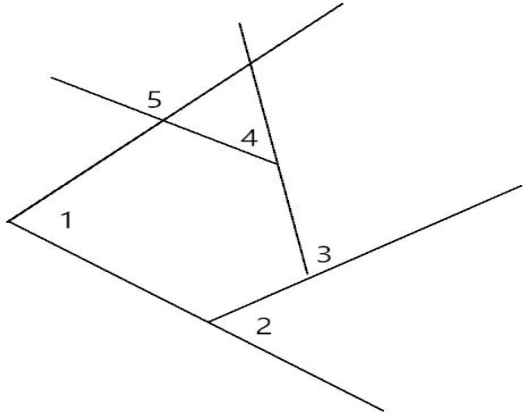

2 Get expressions of every pieces of $B_{jk}^i$ by convolution of cone splines for general B-splines, or by directional integration for Box splines and the corresponding cells.

3. Evaluation of $B_{jk}^i$ on each cell by the corresponding explicit expression.

4 Compute inner product of two polyhedral spline finite elements by summing inner products of each piece of one polyhedral spline finite element with each piece of the other, that is, inner product of two polynomials on a polygon, which can be computed by triangulating the polygon and turning the polynomial into barycentric coordinate forms on each triangle, and then summing the integrals on all the triangles.

If the B-spline is a box spline, we need not turn polygons into combination of cones. We can compute explicit expressions of $B_{jk}^i$ directly. The algorithm is depicted as follows:

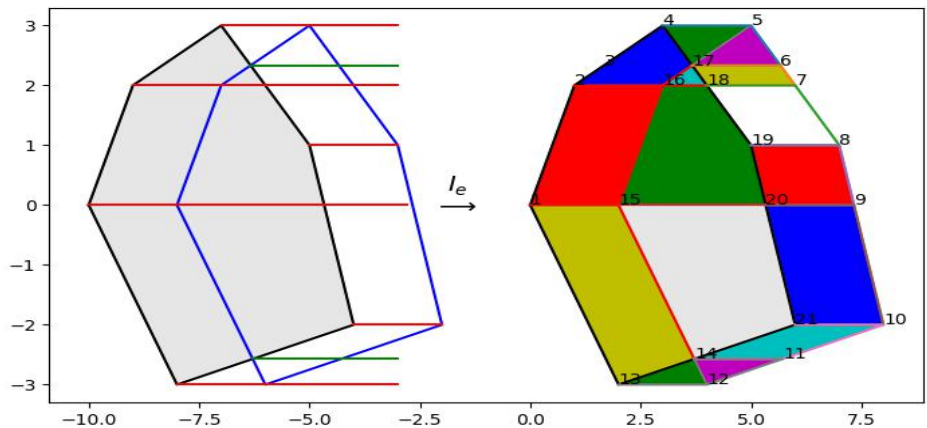

The left big gray polygon is the original cell. Firstly, it is translated by $e=(2,0)$ in the figure to the blue polygon. Secondly, it is sliced horizontally into three cells by the red lines. Thirdly, these three cells are subdivided by the blue lines, and green lines are obtained simultaneously. Lastly, cells outside of the polygon on the right bounded by red, blue and green lines are appended. Cells of the partition after horizontal integration are depicted on the right.

## 4. Subdivision and adaption

Polyhedral spline finite elements are based on arbitrary polygonal partition. The partition can be chosen flexibly. The intersection of two cells may be part of 1-cells. There are more choices for subdivision and adaption. Cells can be subdivided, combined and moved. The degree of $f_{jk}^{i}$ can be raised. The mollifier can be subdivided, degree elevated or replaced. General mesh refinement, moving mesh and adaption are relatively easier to conduct. As is known that it is difficult to construct $c^1$ finite elements, especially for three dimensional cases. There is no such problem for polyhedral spline finite elements. Besides, multigrids and hierachical splines can be used for adaption. These aspects deserve further research, we leave them for future papers.

## 5. Consideration of boundary conditions

If the PDE problem to be solved is a periodic one, the domain can be partitioned periodically. If the PDE problem to be solved is not periodic, in order to apply polyhedral spline finite elements to numerical solution of PDEs, we must enlarge the domain $\Omega$ so that the space of polyhedral spline finite elements keeps approximation power on the boundary of domain $\Omega$. The following figure shows some examples.

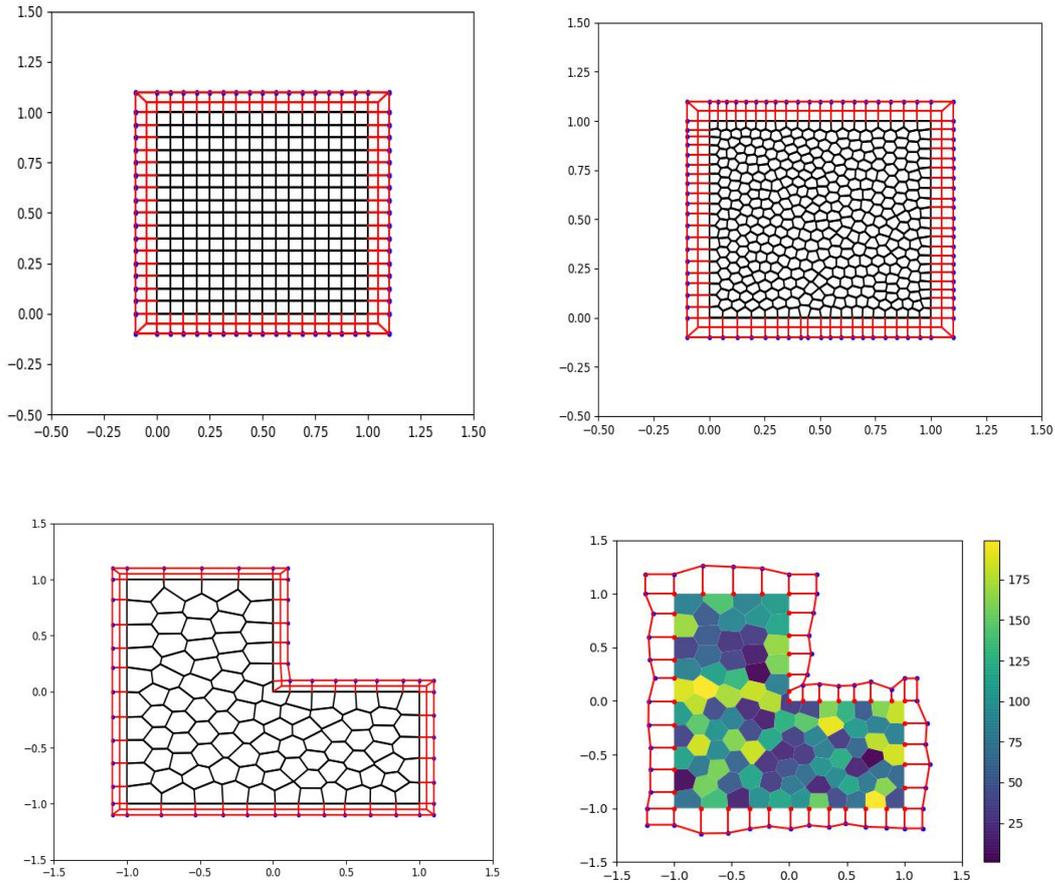

## 6. Combination with isogeometric analysis

Polyhedral spline finite elements and isogeometric analysis can be used together. NURBS are used to transform irregular domains into regular ones. Polyhedral spline finite elements are used for moving mesh and adaption. Polyhedral spline finite elements can be used in place of NURBS as follows

$$L(f) = \sum \lambda_{jk}^{i}(f)\omega_{jk}^{i} B_{jk}^{i} / \sum \omega_{jk}^{i} B_{jk}^{i},$$

since the space of polyhedral spline finite elements contains polynomials of degree at most n, but it is not easy to write out the formula.

## 7. Application in numerical solution of PDEs

According to the discussion of section three, it is not difficult to implement all the algorithms, but it is time consuming, I cannot write out all the programs in one month and a half. So, in this section I give a simple example to demonstrate the application of my method. Now consider the following second order elliptic differential equation:

$$\begin{cases} -\Delta u + \beta \cdot \nabla u + \gamma u = f, & in \quad \Omega \\ \quad\quad u = g, & on \quad \partial\Omega \end{cases}$$

where $\Omega=[0,2n+1]\times[0,2n+1]$. The weak form of the above problem is formulated as follows: find $u \in H^1(\Omega)$ such that $a(u,v) = \int_\Omega fv d\omega$ for all $v \in H_0^1(\Omega)$ where $a(u,v) = \int_\Omega (\nabla u \cdot \nabla v + \beta \cdot \nabla uv + \gamma uv) d\omega$.

Partition $\Omega=[0,2n+1]\times[0,2n+1]$ into blue diamonds and triangles by straight lines $x+y=k, k=1,3,\cdots 4n+1$ and $x-y=k, k=-2n+1,-2n+3,\cdots,2n-1$. Now enlarge $\Omega$ a little and modify the partition according to the following figure. The new partition consists of all the small diamonds: blue ones, red ones and diamonds with both blue and red edges.

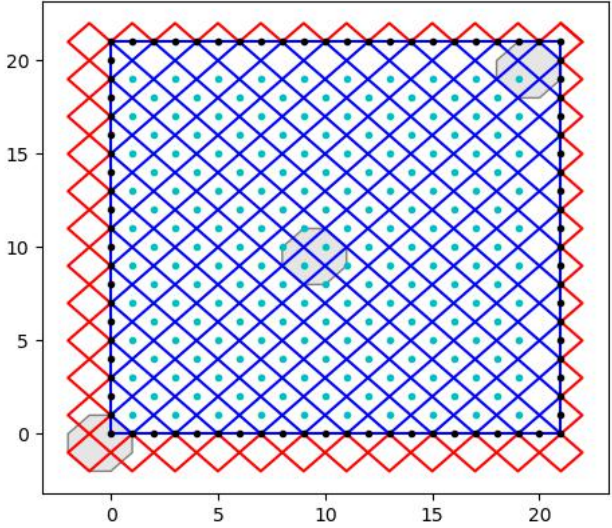

Let $f^i = 1$ on cell $\Delta_i$ and $f^i = 0$ otherwise. Select a B-spline $B(x,y)$, such that $B(x,y)=1$ on $[0,1]\times[0,1]$ and $B(x,y) = 0$ otherwise. Let $B^i = B * f^i$. As is known that $B^i$ is a Box spline in $s_2^1(\Delta_2)$. Now let $u = \sum c_i B^i$ and set $v = B^k$, whose support is contained in $\Omega=[0,2n+1]\times[0,2n+1]$, in the equation $a(u,v) = \int_\Omega fv d\omega$. The diamonds with green dots correspond to Box splines with support in $\Omega$. Set the boundary conditions $u = g$ at black dots. Then a linear system is obtained with $2n^2 + 6n + 5$ unknowns and $2n^2 + 6n + 5$ equations. The concrete computation is not given here. It can be conducted according to section one.

## 8. Summary

Polyhedral spline finite elements have both high smooth order and strong approximation power. They are suitable for moving mesh and adaption. They can be used for geometric modeling. They can be generalized to spline differential forms. Many problems need further study.

**Acknowledgments:** The author would like to express his gratitude to reviewers' comments and remarks that help improve the presentation considerably. This work is supported by the grant from the National Natural Science Foundation of China (No. 61572099).